\documentclass[a4paper,12pt]{article}
\usepackage{amsmath,amsfonts,amssymb,amscd, array}
\usepackage[ps,matrix]{xypic}

\def\N{{\mathbb N}}
\def\Z{{\mathbb Z}}
\def\C{{\mathbb C}}

\def\a{\alpha}
\def\b{\beta}
\def\l{\lambda}

\def\moins{\raise 1pt\hbox{{$\scriptstyle -$}}}
\def\plus{\raise 1pt\hbox{{$\scriptstyle +$}} }

\def\d{\partial}
\def\ph{\widehat{\pi}}
\def\Kh{\widehat K}
\def\Yh{\widehat Y}
\def\Uh{\widehat U}

\def\cH{{\mathcal H}}

\def\mfS{\mathfrak{S}}
\def\p{\mathfrak{p}}
\def\n{\mathfrak{n}}

\def\QED{\hfill $\square$}
\def\Pol{{ \mathfrak{ P\hspace{-0.06 em}o\hspace{-0.05 em} l}}}
\def\Part{{\mathfrak{ P\hspace{-0.06 em}a\hspace{-0.05 em}r\hspace{-0.05 em}t}}}

\def\u{{\mathbf u}}

\newtheorem{example}{Example}[section]

\newtheorem{theorem}[example]{Theorem}
\newtheorem{corollary}[example]{Corollary}
\newtheorem{definition}[example]{Definition}
\newtheorem{proposition}[example]{Proposition}
\newtheorem{lemma}[example]{Lemma}

\def\proof{{\it Proof.}\ }
\def\QED{\hfill Q.E.D}

\begin{document}

 \title{\large\sc Non-Symmetric Hall-Littlewood Polynomials}
\author{Francois Descouens \and Alain Lascoux}
\date{}
\maketitle

\hfill{\it \`A Adriano Garsia, en toute amiti\'e}

\begin{abstract}
Using the action of the Yang-Baxter elements of the Hecke algebra on
 polynomials, we  define two bases of polynomials in $n$ variables. 
The Hall-Littlewood polynomials are a subfamily of one of them.
For $q=0$, these bases specialize into the two families
of classical Key polynomials
(i.e. Demazure characters for type $A$). We  give a scalar product
for which the two bases are adjoint of each other.
\end{abstract}

\section{Introduction}

We define two linear bases of the ring of polynomials in
$x_1,\ldots, x_n$, with coefficients in $q$. 

These polynomials, that we call \emph{$q$-Key polynomials}, 
and denote $U_v,\Uh_v$, $v\in \N^n$,  
specialize at $q=0$ into key polynomials $K_v, \Kh_v$. 
The polynomials $U_v$ which are symmetrical in $x_1,\ldots, x_n$
are precisely the Hall-Littlewood polynomials $P_\l$, 
indexed by partitions $\l\in \Part$, the relation between the two indices
being  $\l=[\l_1,\ldots,\l_n]= [v_n,\ldots,v_1]$.

Our main tool is the  Hecke algebra $\cH_n(q)$ 
of the symmetric group,
acting on polynomials by deformation of divided differences.
This algebra  contains two adjoint bases of Yang-Baxter 
elements (Th. \ref{Flag}). 
The $q$-Key polynomials are the images of dominant monomials under
these Yang-Baxter elements (Def. \ref{DefKey}).
These polynomials are clearly two linear bases of polynomials,
since the transition matrix to monomials is uni-triangular.

We show in the last section that $\{U_v\}$ and $\{\Uh_v\}$ 
are two adjoint bases with respect to a certain scalar product
reminiscent of Weyl's scalar product on symmetric functions.

We have intensively used MuPAD (package {\tt MuPAD-Combinat} \cite{MuPAD})
and Maple  (package {\tt ACE} \cite{ACE}).

\section{The Hecke algebra $\cH_n(q)$}

Let $\cH_n(q)$ be the Hecke algebra of the symmetric group $\mfS_n$,
with coefficients the rational functions in a parameter $q$.
It has
generators $T_1,\ldots,T_{n-1}$ satisfying the braid relations 
\begin{equation}
\begin{cases}
 T_iT_{i+1}T_i   =  T_{i+1}T_iT_{i+1},  \\
 T_iT_j=T_jT_i ~~ (\vert j-i \vert > 1)\, ,
\end{cases}
\end{equation}
and the Hecke relations
\begin{equation}
(T_i+1)(T_i-q)=0 \ ,\ 1\leq i\leq n-1 
\end{equation}

For a permutation $\sigma$ in $\mfS_n$, we denote by $T_\sigma$ the
element $ T_\sigma=T_{i_1}\ldots T_{i_p}$ where $(i_1,\ldots,i_p)$ is
any reduced decomposition of $\sigma$. 
The set $\{T_\sigma:\, \sigma\in\mfS_n\}$ is a linear basis of
$\cH_n(q)$.

\subsection{Yang-Baxter bases}

Let $s_1,\ldots, s_{n-1}$ denote the simple transpositions,
$\ell(\sigma)$ denote the length of $\sigma\in\mfS_n$,
and let $\omega$ be the permutation of maximal length.

Given any set of indeterminates
$\u=(u_1,\ldots,u_n)$, let
$\cH_n(q)[u_1,\ldots,u_n] = \cH_n(q) \otimes \C[u_1,\ldots,
u_n]$. 

One defines recursively 
a \emph{Yang-Baxter basis}  $(Y^{\u}_\sigma)_{\sigma
\in \mfS_n}$, depending on $\u$,  by
\begin{equation}
Y^\u_{\sigma s_i} = Y^\u_{\sigma}\, \left(T_i+ \frac{1-q}{ 1-
u_{\sigma_{i+1}}/ u_{\sigma_i}} \right), \quad \text{when}\ \ell(\sigma
s_i) > l(\sigma) \, ,
\end{equation}
starting with $Y_{id}^\u = 1$.

Let $\varphi $ be the anti-automorphism of $\cH_n(q)[u_1,\ldots,
u_n]$ such that
$$ \left \lbrace
\begin{array}{l}
\varphi(T_\sigma)=T_{\sigma^{-1}},\\
\varphi(u_i)=u_{n-i+1}.
\end{array} \right .$$

We define a bilinear form $<\, ,\, >$ 
on $\cH_n(q)[u_1,\ldots,u_n]$ by
\begin{equation}
<h_1\, ,\, h_2>\ :=\ 
\text{coefficient of }T_\omega\ \text{in}\ 
h_1\cdot \varphi(h_2) \ .
\end{equation}

The main result of \cite[Th. 5.1]{FlagYang} 
is the following duality property of Yang-Baxter bases.
\begin{theorem}  \label{Flag}
For any set of parameters $\u=(u_1,\ldots,u_n)$, the basis adjoint to
$(Y^{\u}_\sigma)_{\sigma\in \mfS_n}$ with respect to $<\, ,\, >$ is
the basis $(\Yh^\u_\sigma)_{\sigma \in
\mfS_n}=(Y^{\varphi(\u)}_\sigma)_{\sigma\in \mfS_n}$. More precisely,
one has
$$ \forall \ \sigma\,, \nu\, \in \mfS_n\,, \quad <Y^\u_\sigma\, , \,
\Yh^\u_\nu>=\delta_{\lambda,\ \nu\omega} \quad .$$
\end{theorem}

Let us fix from now on the parameters $u$ to be $ \u = (1,q,q^2,\ldots,
q^{n-1})$. Write $\cH_n$ for 
$\cH_n(q)[1,q,\ldots,q^{n-1}]$.

In that case, the Yang-Baxter basis
$(Y_{\sigma})_{\sigma\in \mfS_n}$ and its adjoint basis
$(\Yh_{\sigma})_{\sigma\in \mfS_n}$ are defined recursively,
starting with $Y_{id}=1=\Yh_{id}$, 
 by
\begin{equation}  \label{YBdef}
 Y_{\sigma s_i} = Y_{\sigma}\, \left(T_i + 1/\lbrack k
\rbrack_q \right) \quad \text{and}\,\, \Yh_{\sigma s_i} =
\Yh_{\sigma}\, \left(T_i + q^{k-1}/ \lbrack k \rbrack_q
\right)\,)\, ,\ \ell(\sigma s_i)> \ell(\sigma)    \,,
\end{equation}
with  $k=\sigma_{i+1} -\sigma_i$ and 
$\lbrack k\rbrack_q= (1-q^k)/(1-q)$.

Notice that  the maximal Yang-Baxter elements 
have another expression \cite{DKLLST}~:
$$Y_\omega = \sum_{\sigma\in \mfS_n}\, T_\sigma \ \ \text{and} \ \
\Yh_\omega = \sum_{\sigma\in \mfS_n} \, (-q)^{\ell(\sigma\omega)}\,
T_\sigma \ . $$

\begin{example} 
For $\cH_3$, the transition matrix 
between $\{ Y_\sigma\}_{\sigma\in \mfS_3}$   and 
$\{ T_\sigma\}_{\sigma \in \mfS_3}$  is
$$ \begin{array}{c|cccccc|}
       \text{\footnotesize $123$ \normalsize} & 1     & 1     & 1     & \frac{1}{q + 1} & \frac{1}{q + 1} & 1 \\
       \text{\footnotesize $132$ \normalsize} &\cdot & 1     & \cdot & 1               & \frac{1}{q + 1} & 1   \\
       \text{\footnotesize $213$ \normalsize} &\cdot & \cdot & 1     & \frac{1}{q + 1} & 1               & 1 \\
       \text{\footnotesize $231$ \normalsize} &\cdot & \cdot & \cdot & 1               & \cdot           & 1 \\
       \text{\footnotesize $312$ \normalsize} &\cdot & \cdot & \cdot & \cdot           &          1      & 1 \\
       \text{\footnotesize $321$ \normalsize} &\cdot & \cdot & \cdot 
      & \cdot    & \cdot  & 1 \end{array} \quad ,$$
writing $`\cdot`$ for $0$. Each column represents the expansion
of some element  $Y_\sigma$.
\end{example}

\subsection{Action of $\cH_n$ on polynomials}

Let $\Pol$ be the ring of polynomials in the variables
$x_1,\ldots,x_n$ with coefficients the rational functions in
$q$. We write monomials exponentially: $x^v= x_{1}^{v_1}\ldots x_{n}^{v_n}$,
 $v=(v_1,\ldots,v_n) \in \mathbb{Z}^n$. 
A monomial $x^v$ is \emph{dominant} if
$v_1 \ge \ldots \ge v_n$. 

We extend the natural order on partitions to elements of $\mathbb{Z}^n$
by
$$u\le v\, \quad \text{iff}\quad \forall k> 0\,,\quad
\sum_{i=k}^{n}(v_{i}-u_{i}) \ge 0\, .$$ 
For any polynomial
$P$ in $\Pol$, we call {\it leading term} of $P$ all the 
monomials (multiplied by their coefficients) 
which are maximal with respect to this partial order.
This order is compatible with
 the right-to-left lexicographic order, that we shall also use. 
We also use the classiccal notation 
$\n(v) =0v_1+1v_2+2v_3+\cdots +(n\moins 1)v_n$.

Let $i$ be an integer such that $1\le i \le n-1$. As an operator on
$\Pol$, the simple transposition $s_i$ acts by 
switching $x_i$ and $x_{i+1}$, and we denote this action by 
$f\to f^{s_i}$.
The $i$-th
\emph{divided difference} $\d_i$ and the $i$-th 
\emph{isobaric divided difference} 
$\pi_i$, written on the right of the operand,  are the 
following operators~:
$$\d_i: \ f \longmapsto f\, \d_i :=\frac{f -f^{s_i}}{x_i-x_{i+1}} \qquad  ,
\qquad \pi_i: \ f \longmapsto  f\, \pi_i :=\frac{x_if
-x_{i+1}f^{s_i}}{x_i-x_{i+1}} \ . $$

 The Hecke algebra $\cH_n$ has a faithful representation  as an algebra of
operators on $\Pol$ 
given by the following equivalent formulas \cite{DKLLST, Lusztig}
$$\left \lbrace \begin{array}{lllllll} T_i &=& \square_i-1 &=&
(x_i-qx_{i+1})\, \d_i - 1 &=& (1-qx_{i+1}/x_i)\pi_i-1\, ,\\ Y_{s_i}
&=&\square_i &=& (x_i-qx_{i+1})\, \d_i &=& (1-qx_{i+1}/x_i)\pi_i \, ,
\\ \hat{Y}_{s_i} &=&\nabla_i &=& \square_i - (1+q) &=& \d_i\,
(x_{i+1}-qx_i) \, .
\end{array}\right .$$

The Hecke relations  imply
$$ \square_i^2=(1+q)\square_i ~~~,~~~ \nabla_i^2=-(1+q)\nabla_i
 ~~~\text{and}~~~ \square_i\nabla_i=\nabla_i\square_i=0 \ .$$ 
One easily checks that 
the operators $R_i(a,b)$ and $S_i(a,b)$ defined by
$$ R_i(a,b)=\square_i-q\frac{[b-a-1]_q}{[b-a]_q}\quad \text{and} \quad
S_i(a,b)=\nabla_i+q\frac{[b-a-1]_q}{[b-a]_q}$$ satisfy the Yang-Baxter
equation
\begin{equation}
R_i(a,b)\ R_{i+1}(a,c) \ R_i(b,c) = R_{i+1}(c,b)\ R_i(a,c) \
R_{i+1}(a,b) \, .
\end{equation}
We have implicitely used these equations in the recursive definition
of Yang-Baxter elements (\ref{YBdef}). 

This realization comes from geometry \cite{chi_y},
where the maximal Yang-Baxter
elements are interpreted as Euler-Poincar\'e characteristic for the
flag variety of $GL_n(\C)$.  
This gives still another expression of the maximal Yang-Baxter
elements~:
\begin{equation}  \label{MaxYang}
 Y_\omega = \prod_{1\leq i<j\leq n} (x_i-q x_j)\, \d_\omega 
\qquad ,\qquad 
\Yh_\omega = \d_\omega\, \prod_{1\leq i<j\leq n} (x_j-q x_i) \, .
\end{equation}

\begin{example} Let $\sigma=(3412)=s_2s_3s_1s_2$.
The elements $Y_{3412}$ and $\Yh_{3412}$ can be written
\begin{eqnarray*}
Y_{3412}&=&\square_2\left(\square_3-\frac{q}{1+q}\right)\left(\square_1-\frac{q}{1+q}\right)
\left(\square_2-\frac{q+q^2}{1+q+q^2}\right)\ ,\\
\Yh_{3412}&=&\nabla_2\left(\nabla_3+\frac{q}{1+q}\right)\left(\nabla_1+\frac{q}{1+q}\right)
\left(\nabla_2+\frac{q+q^2}{1+q+q^2}\right)\ .
\end{eqnarray*}
\end{example}

We shall now identify the images of dominant monomials under the
maximal Yang-Baxter operators with Hall-Littlewood polynomials.
Recall that there are two proportional families $\{ P_\l \}$ 
and $\{ Q_\l \}$ of Hall-Littlewood polynomials.
Given a partition 
$\lambda=[\l_1,\l_2,\ldots, \l_r]= 
(0^{m_0},1^{m_1},\ldots,n^{m_n})$, with $m_0=n-r=n-m_1-\cdots -m_n$,
then  
$$ Q_\l = \prod_{1\leq i \leq n} \prod_{j=1}^{m_i} (1-q^j)\, P_\l \, .$$

Let moreover $d_\l(q)= \prod_{0\leq i \leq n} \prod_{j=1}^{m_i}[j]_q$. 
The definition of Hall-Littlewood polynomials with raising operators
\cite{Littlewood},\cite[III.2]{Macdonald} 
can be rewritten, thanks to (\ref{MaxYang}),
as follows.

\begin{proposition}
Let $\lambda$ be a partition of $n$.
Then one has
\begin{equation}
x^\lambda \ Y_\omega\, d_\lambda(q)^{-1} = 
P_{\lambda}(x_1, \ldots, x_n;q)
\end{equation}
\end{proposition}

The family of the Hall-Littlewood functions $\{Q_\l \}$
indexed by  partitions can be  extended  into a family 
$\{ Q_v :\, v\in \mathbb{Z}^n\}$, 
using the following relations due to Littlewood (\cite{Littlewood},
\cite[II.2.Ex. 2]{Macdonald})
\begin{equation}\label{rule1}
Q_{(\ldots,u_{i},u_{i+1},\ldots)}=
-Q_{(\ldots,u_{i+1}-1,u_{i}+1,\ldots)}
+q\ Q_{(\ldots,u_{i+1},u_{i},\ldots)}
+q\ Q_{(\ldots,u_{i}+1,u_{i+1}-1,\ldots)}\quad\text{if}\quad u_i<u_{i+1} ,\\ 
\end{equation}
\begin{equation}\label{rule2}
Q_{(u_1,\ldots,u_n)}=0 \quad \text{if}\quad u_n<0 \quad .
\end{equation}
By iteration of the first relation, one can write any $ Q_u$ in terms
of Hall-Littlewood functions indexed by decreasing vectors $v$ such
that $\vert v \vert=\vert u\vert$. Consequently, if $u$ such
that $\vert u \vert=0$, $Q_u$ must be proportionnal to 
$Q_{0\ldots 0}=1$, i.e. is a constant that one can note as the
specialisation $Q_u(0)$ in $x_1=0=\cdots=x_n$.

The final expansion of $Q_u$, after iterating (\ref{rule1}) many times,
is  not easy to predict. In particular, one needs to know whether
$Q_u\neq 0$.  For that purpose, we shall isolate a distinguished
term in the expansion of $Q_u$. Given a sum 
$\sum_{\l\in \Part} c_\l(t) Q_\l$, call \emph{top term}
the image of the leading term $\sum c_\mu(t) Q_\mu$ after 
restricting each coefficient $c_\mu(t)$ 
to its term in highest degree in $t$. 

Given $u\in\Z^n$, define recursively $\p(u)\in \Part \cup \{-\infty\}$ by
\begin{itemize}
\item if $u\not\geq [0,\ldots,0]$ then $\p(u)=-\infty$
\item if $u_2\geq u_3\geq \cdots \geq u_n>0$ then
  $\p(u)$ is the maximal partition of length $\leq n$, of weight
  $|u|$ (eventual zero terminal parts are suppressed).
\item $\p(u)= \p\bigl(u\, \p([u_2,\ldots, u_n])    \bigr)$  
\end{itemize}

\begin{lemma}   \label{NonNul2}
Let $u\in\mathbb{Z}^n$. Then
\begin{itemize}
\item if $u\not\geq [0,\ldots,0]$ then $Q_u=0$\, ,
\item if $u\geq [0,\ldots,0]$, let $v=\p(u)$. 
 Then $Q_u\neq 0$ and its leading term is $q^{\n(u)-\n(v)} Q_v$\, .
\end{itemize}

\end{lemma}

\proof    Given any  decomposition $u=u'.u''$, then one can 
apply (\ref{rule1}) to $u''$ and write $Q_u$ as a linear combination
of terms $Q_{u' v}$ with $v$ decreasing, with $|v|=|u''|$. 
Therefore, if $|u''|=0$, then the last components of such $v$
are negative, all $Q_{u' v}$ are $0$, and $Q_u=0$.

If $u \geq [0,\ldots,0]$ and $u$
 is not a partition, write $u=[\ldots, a,b,\ldots]$,
with $a,b$ the rightmost increase in $u$. 
We apply relation (\ref{rule1}), assuming the validity of lemma 
for the three terms in the RHS ~:
$$  Q_{\ldots, a,b,\ldots} =
  -Q_{\ldots, b-1,a+1,\ldots}
+q  Q_{\ldots, b,a,\ldots}+  q Q_{\ldots, a+1,b-1,\ldots}   $$
Notice that the first two terms have not necessarily an index
$\geq [0,\ldots,0]$, but that $[\ldots, a+1,b-1,\ldots]
\geq [0,\ldots,0]$. 

In any case, it is clear that 
$\p([\ldots, b-1,a+1,\ldots])=p_1\leq v$, 
$\p([\ldots, b,a,\ldots])=p_2\leq v$,
and $\p([\ldots, a+1,b-1,\ldots])=v$.  

Restricted to top terms, the expansion of the RHS in the basis $Q_\l$
becomes
$$ - \Bigl( (q^{\n(u)+a+1-b-\n(v)}+\cdots)\, Q_v   \Bigr) 
   + q \Bigl( (q^{\n(u)+a-b-\n(v)}+\cdots)\, Q_v   \Bigr)
   +q \Bigl( (q^{\n(u) -1-\n(v)}+\cdots)\, Q_v   \Bigr) \, ,$$
where one or two of the first two terms may be replaced by $0$,
depending on the value of $p_1$, or $p_2$.
In final, the top term of the RHS is 
$ q^{\n(u)-\n(v)}Q_v$, as wanted.                  \QED

\begin{example} For $v=[-2,3,2]$,
$$ Q_{2,3,2} = (q^3-q^2) Q_3 + (q^5+q^4-q^3-2q^2+q) Q_{21}
  + (q^4-q^3-q^2+q) Q_{111}  \, ,$$
and the top term is $q^4 Q_{111}$,
since $4=(0(-2)+1(3)+2(2)) -(0(1)+1(1)+2(1))$ and $[1,1,1]>[2,1]$,
$[1,1,1]>[3]$. Notice that the coefficient of $Q_{21}$ is of 
higher degree.
\end{example}

\section{$q$-Key Polynomials}

In this section, 
we show that the images of dominant monomials under the 
Yang-Baxter elements $Y_\sigma$ (resp. $\Yh_\sigma$),
$\sigma\in\mfS_n$ constitute two bases of $\Pol$, 
which specialize into the two families of Demazure characters.

We have already identified in the preceding section  the images 
of dominant monomials under $Y_\omega$ 
to Hall-Littlewood polynomial, using the relation between 
 $Y_\omega$ and $\d_\omega$. The other polynomials are new.  

\subsection{Two bases}

The dimension of the linear span of the image of a monomial $x^v$
 under all permutations depends upon the stabilizer of $v$.
We meet the same  phenomenon  when taking 
the images of a monomial under Yang-Baxter elements.

Let $\l=[\l_1,\ldots, \l_n]$ be a decreasing partition
(adding eventual parts equal to $0$).
Denote its orbit under permutations of components  by $\mathcal{O}(\lambda)$.
Given any $v$ in $\mathcal{O}(\lambda)$, let $\zeta(v)$ be the
permutation of maximal length such that $\l\, \zeta(v) = v$ and
$\eta(v)$ be the permutation of minimal length such that $\l\, \eta(v)
= v$.  These two permutations are representative of the same coset of
$\mfS_n$ modulo the stabilizer of $\lambda$.
\begin{definition}  \label{DefKey}
For all $v$ in $\N^n$, the \emph{$q$-Key polynomials} $U_v$ and
$\Uh_v$  are the following polynomials~:
$$ U_v(x;q)= \left (\frac{1}{d_\l(q)} x^\l\right ) Y_{\zeta(v)} \qquad
, \qquad \Uh_v(x;q)= x^\l \Yh_{\eta(v)} \, ,$$ where $\lambda$
is the dominant reordering of $v$.
\end{definition}
In particular, if $v$ is (weakly) increasing, then 
$\zeta(v)=\omega$ and  $U_v$
is a Hall-Littlewood polynomial.

\begin{lemma}  \label{Transition2} 
The leading term of $U_v$ and  $\Uh_v$ is $x^v$.
Consequently, 
the transition matrix between the $U_v$ (resp. the $\Uh_v$) and the
monomials is upper unitriangular with respect to the 
right-to-left lexicographic order.
\end{lemma}

\proof    
  Let $k$ be an integer and $u$ be a weight such that
$u_k>u_{k+1}$. Suppose by induction that
$x^u$ is the leading term of $U_u$.  
Recall the the explicit action of $\square_k$ is 
(noting only the two variables $x_k,x_{k+1}$)        
\begin{eqnarray*} 
   x^{\b\a}\, \square_k&=& x^{\b\a} +(1-t)( x^{\b-1,\a+1}+\cdots +
                  x^{\a+1,\b-1}) + x^{\a\b} \, , \ \b>\a    \\
   x^{\b\b}\, \square_k&=&  (1+t) x^{\b\b}   \\ 
 x^{\a\b}\, \square_k&=& tx^{\b\a} +(t-1)( x^{\b-1,\a+1}+\cdots +
                  x^{\a+1,\b-1}) + tx^{\a\b} \, , \ \a<\b \, . 
\end{eqnarray*}
>From  these formulas, it is clear that
for any constant $c$, the leading term of 
$x^u\, (\square_k+c)$  is $(x^u)^{s_k}$, 
and, for any $v$ such that $v<u$, all the monomials 
in $x^v\, (\square_k+c)$ are strictly less (with respect to the partial
order) than $(x^u)^{s_k}$.  \hfill $\square$

\begin{example}
For $n=3$, Figures {\emph 1}
and {\emph 2} 
 show the  case of a regular dominant weight $x^{210}$ and
Figures {\emph 3} and {\emph 4} correspond to a case, $x^{200}$, 
where the stabilizer is not trivial. 
In this last case, the polynomials belonging to the 
family are framed, the extra polynomials denoted $A,B$ 
do not belong to the basis. 
\begin{figure}[!h]
\newdimen\vcadre\vcadre=0.1cm 
\newdimen\hcadre\hcadre=0.1cm 
\def\GrTeXBox#1{\vbox{\vskip\vcadre\hbox{\hskip\hcadre\scriptsize
                       $#1$
                      \hskip\hcadre}\vskip\vcadre}}
$\vcenter{\xymatrix@R=2cm@C=-3cm{ & * {\GrTeXBox{U_{210}=x^{210}}}
\ar @{->}[dl]^{\square_1} \ar @{->}[dr]^{\square_2}& \\
*{\GrTeXBox{U_{120}=x^{120}+x^{210}}} \ar @{->}[d]^{\square_2-q/(1+q)}
& & *{\GrTeXBox{U_{201}=x^{201}+x^{210}}} \ar
@{->}[d]^{\square_1-q/(1+q)} \\
*{\GrTeXBox{U_{102}=x^{102}+(1-q)x^{111}+\frac{1}{1+q}x^{120}+x^{201}+\frac{1}{1+q}x^{210}}
} \ar @{->}[dr]^{\square_1} & & *{\GrTeXBox{U_{021}=x^{021}+(1-q)x^{111}+\frac{1}{1+q}x^{201}+x^{120}+\frac{1}{1+q}x^{210}}} \ar @{->}[dl]^{\square_2} \\ & 
*{\GrTeXBox{U_{012}= 
+ x^{102} + (- q - q + 2)x^{111} + x^{120} + x^{201} + x^{210} } } & \\ }} $
\caption{$q$-Key polynomials generated from $x^{210}$.}
\end{figure}
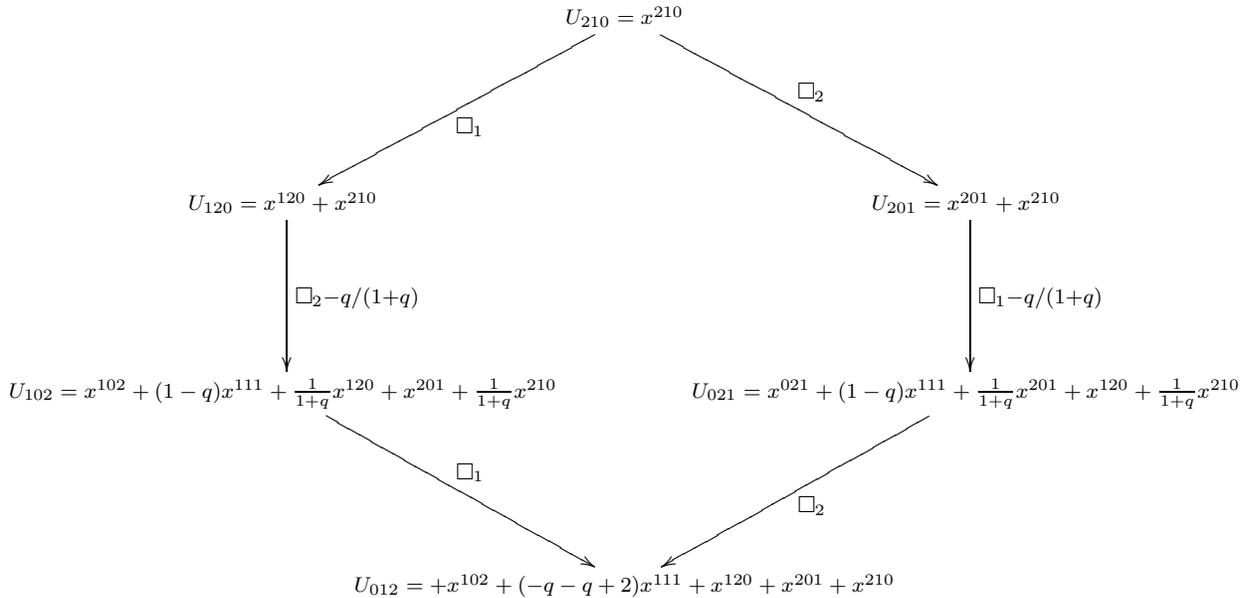
\begin{figure}[!h]
\newdimen\vcadre\vcadre=0.1cm 
\newdimen\hcadre\hcadre=0.1cm 
\def\GrTeXBox#1{\vbox{\vskip\vcadre\hbox{\hskip\hcadre \scriptsize
                       $#1$
                      \hskip\hcadre}\vskip\vcadre}}
$\vcenter{\xymatrix@R=2cm@C=-1cm{ & *
{\GrTeXBox{\Uh_{210}=x^{210}}} \ar @{->}[dl]^{\nabla_1} \ar
@{->}[dr]^{\nabla_2}& \\ *{\GrTeXBox{\Uh_{120}=x^{120}-qx^{210}}} \ar
@{->}[d]^{\nabla_2+\frac{q}{1+q}}& &
*{\GrTeXBox{\Uh_{201}=x^{201}-qx^{210}}} \ar
@{->}[d]^{\nabla_1+\frac{q}{1+q}} \\
*{\GrTeXBox{\Uh_{102}=x^{102}+(1-q)x^{111}-qx^{201}-\frac{q^2}{q+1}x^{120}+\frac{q^3}{q+1}x^{210}}}
\ar @{->}[dr]^{\nabla_1} & & *{\GrTeXBox{
\Uh_{021}=x^{021}+(1-q)x^{111}-qx^{120}-\frac{q^2}{q+1}x^{201}+\frac{q^3}{q+1}x^{210}}}
\ar @{->}[dl]^{\nabla_2} \\ & *{\GrTeXBox{\Uh_{012}= (1+q)
B-q\Kh_{120}} } & \\ }} $
\caption{Dual $q$-Key polynomials generated from $x^{210}$.}
\end{figure}
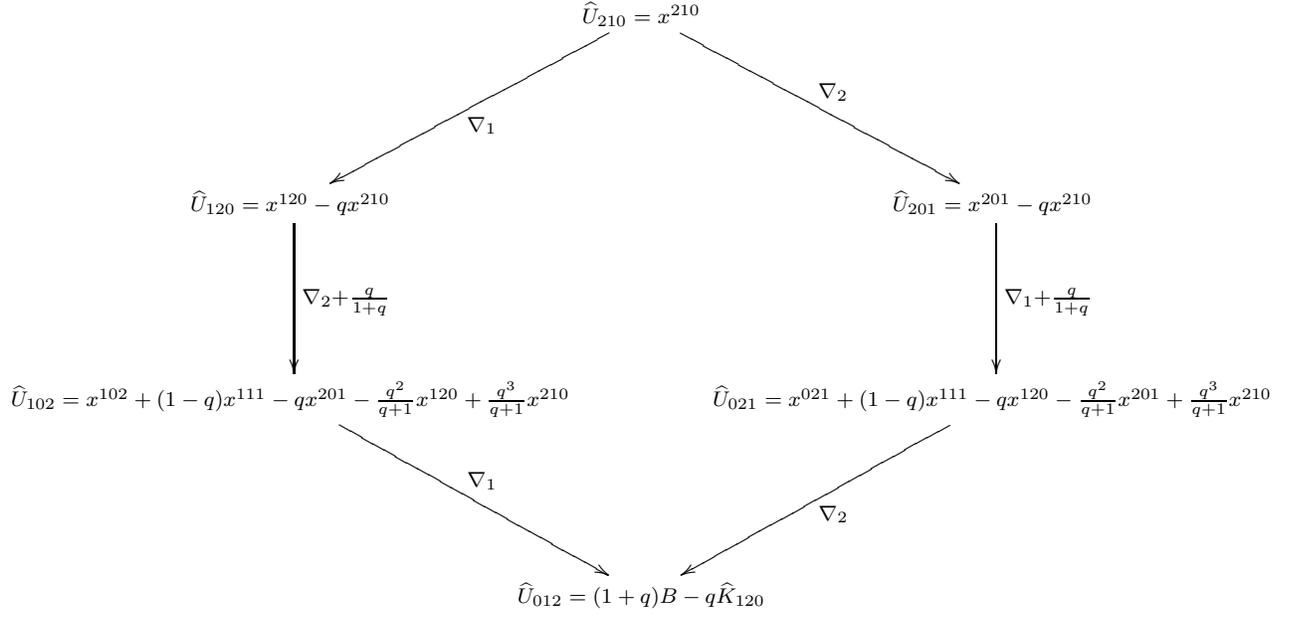

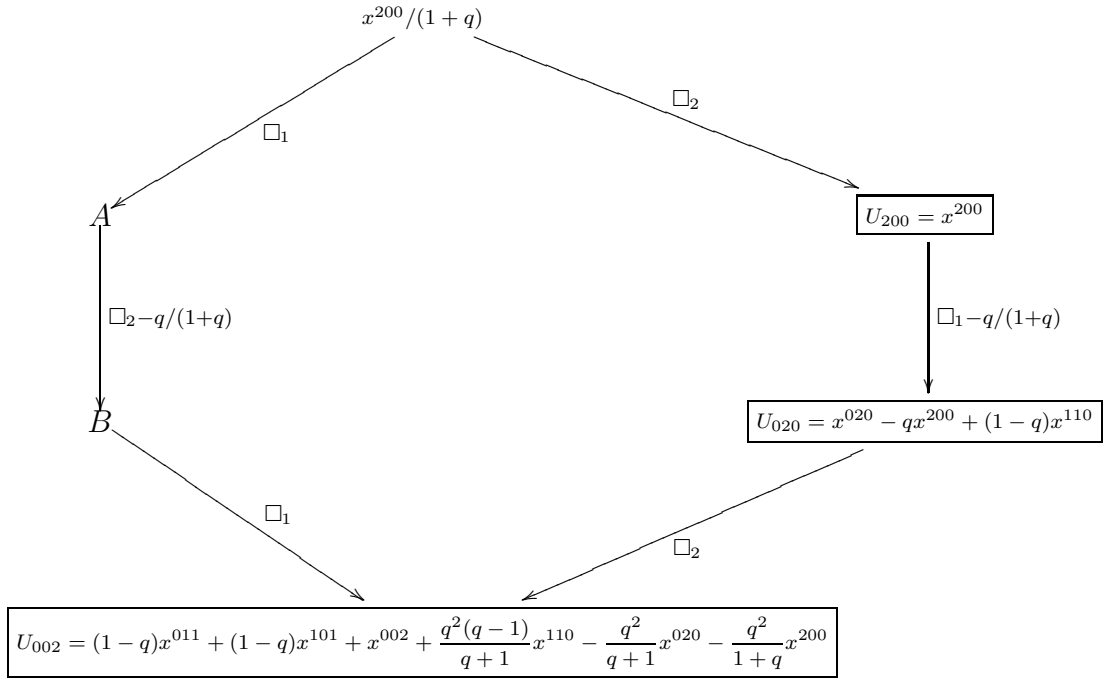
\begin{figure}[!ht]
\newdimen\vcadre\vcadre=0.1cm 
\newdimen\hcadre\hcadre=0.1cm 
\def\GrTeXBox#1{\vbox{\vskip\vcadre\hbox{\hskip\hcadre \scriptsize
                       $#1$
                      \hskip\hcadre}\vskip\vcadre}}
$\vcenter{\xymatrix@R=2cm@C=-1.5cm{ & * 
{\GrTeXBox{x^{200}/(1+q)}}
  \ar @{->}[dl]^{\square_1} 
  \ar @{->}[dr]^{\square_2}& \\ *{A} 
  \ar @{->}[d]^{\square_2-q/(1+q)} & & *
{\GrTeXBox{\boxed{U_{200} =x^{200}}}} 
  \ar @{->}[d]^{\square_1-q/(1+q)} \\*
{B} 
  \ar @{->}[dr]^{\square_1} & &*
{\GrTeXBox{\boxed{U_{020}=x^{020}-qx^{200}+(1-q)x^{110}}}} 
\ar   @{->}[dl]^{\square_2} \\ & *
{\GrTeXBox{\boxed{U_{002}=(1-q)x^{011}+(1-q)x^{101}+x^{002}+\frac{q^2(q-1)}{q+1}x^{110}-\frac{q^2}{q + 1}x^{020}-
\frac{q^2}{1+q}x^{200} }}}
 & \\ }}$ 
\caption{$q$-Key polynomials generated from $x^{200}/(1+q)$.}
\end{figure}
\begin{figure}[!h]
\newdimen\vcadre\vcadre=0.1cm 
\newdimen\hcadre\hcadre=0.1cm 
\def\GrTeXBox#1{\vbox{\vskip\vcadre\hbox{\hskip\hcadre \scriptsize
      $#1$ \hskip\hcadre}\vskip\vcadre}}
$\vcenter{\xymatrix@R=1.3cm@C=0.3cm{&
    *{\GrTeXBox{\boxed{\Uh_{200}
    =x^{200}}}} \ar @{->}[dl]^{\nabla_1} \ar @{->}[dr]^{\nabla_2} \\
    *{\GrTeXBox{\boxed{\Uh_{020}
   =(1-q)x^{110}+x^{020}+\frac{1}{q+1}x^{200} }}} \ar
    @{->}[d]^{\nabla_2+q/(1+q)} & & {\GrTeXBox{0}}\ar @{->}[d]^{\nabla_1+q/(1+q)} \\
    *{\GrTeXBox{\boxed{\Uh_{002}
  =(1-q)x^{011}+(1-q)x^{101}+x^{002}+(1-q)x^{110}+x^{020}+x^{200}}}}
        \ar @{->}[dr]^{\nabla_1}  & &  {\GrTeXBox{0}}\ar @{->}[dl]^{\nabla_2} \\ &* {\GrTeXBox{0}}&  } } $
\caption{Dual $q$-Key polynomials generated from $x^{200}$.}
\end{figure}
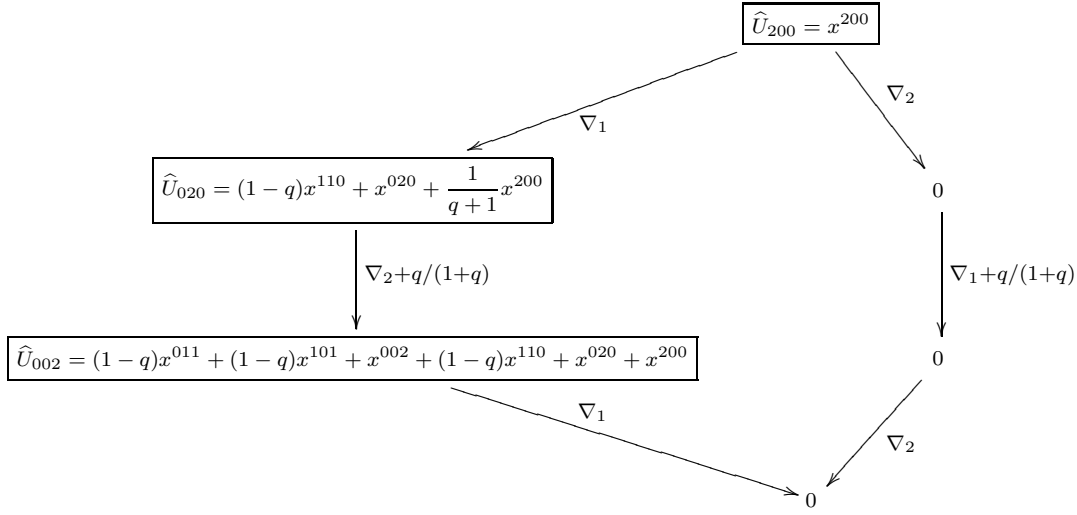
\normalsize
\end{example}

\subsection{Specialization at $q=0$}
The specialization at $q=0$ of the Hecke algebra is called the
  \emph{$0$-Hecke algebra}.
The elementary Yang-Baxter elements specialize in that case into 
\begin{eqnarray}
 Y_{s_i}=T_i+1={\square_i} &\to& x_i\d_i =\pi_i\ , \\
 \hat{Y}_{s_i}=T_i={\nabla_i} &\to&  \d_i x_{i+1} =\ph_i .
\end{eqnarray}

\begin{definition}[Key polynomials]
Let $v\in \N^n$. 
The \emph{Key polynomials} $K_v$ and
$\widehat{K}_v$ are defined recursively, 
starting with $K_v=x^v =\widehat{K}_v$ if $x^v$ dominant, 
 by 
$$ K_{v s_i} = K_v\, \pi_i\qquad , \qquad 
 \widehat{K}_{v s_i} = \Kh_v\, \ph_i\ ,\quad  
\text{for $i$ such that } v_i>v_{i+1}\, .$$
\end{definition}

In particular, the subfamily $(K_v)$ for $v$ increasing, is the
family of Schur functions in $x_1,\ldots,x_n$. 
Demazure \cite{Demazure} defined Key polynomials (using another terminology)
for all the classical groups, and not only the type $A_{n-1}$
which is our case.

Lemma \ref{Transition2} specializes into~:

\begin{lemma}   \label{Transition3}
The transition matrix between the $U_v$ and the $K_v$
(resp. from $\Uh_v$ to $\Kh_v$) is upper unitriangular with
respect to the lexicographic order.
\end{lemma}

\begin{example}
For $n=3$,  the transition matrix
between $\{U_v\}$ and $\{K_v\}$ in weight 3 is (reading a column
as the expansion of some $U_v$) \footnotesize
$$ \begin{array}{c|cccccccccc|} 
\text{\footnotesize $300$ \normalsize}&1 & \cdot & \cdot & \cdot & \cdot & \cdot &\frac{-q}{\left(q + 1\right)} & \cdot & \cdot & \cdot\\
\text{\footnotesize $210$ \normalsize}&\cdot & 1 & \cdot & \cdot & \cdot & \cdot & \cdot & \cdot & \cdot & \cdot\\ 
\text{\footnotesize $201$ \normalsize}&\cdot & \cdot & 1 & \cdot & \cdot & \cdot & \cdot & \frac{-q}{\left(q +1\right)} & \cdot & \cdot\\ 
\text{\footnotesize $120$ \normalsize}& \cdot & \cdot & \cdot & 1 & \cdot & \frac{-q}{\left(q + 1\right)} & -q & \cdot & \cdot & \cdot\\ 
\text{\footnotesize $111$ \normalsize}&\cdot & \cdot & \cdot & \cdot & 1 & -q & \cdot & -q & -q\left(q + 1\right) & q^2\\
\text{\footnotesize $102$ \normalsize}&\cdot & \cdot & \cdot & \cdot & \cdot & 1 & \cdot & \cdot & \cdot & \cdot\\ 
\text{\footnotesize $030$ \normalsize}&\cdot & \cdot & \cdot & \cdot & \cdot & \cdot & 1 & \cdot & \cdot & \cdot\\ 
\text{\footnotesize $021$ \normalsize}&\cdot & \cdot & \cdot & \cdot & \cdot & \cdot & \cdot & 1 & \cdot & \cdot\\ 
\text{\footnotesize $012$ \normalsize}&\cdot & \cdot & \cdot & \cdot & \cdot & \cdot & \cdot & \cdot & 1 & -q\\ 
\text{\footnotesize $003$ \normalsize}&\cdot & \cdot & \cdot & \cdot & \cdot & \cdot & \cdot & \cdot & \cdot & 1\end{array}\, ,$$ \normalsize and the
 transition matrix between  $\{\Uh_v\}$ and  $\{\Kh_v\}$ is \footnotesize
$$ 
\begin{array}{c|cccccccccc|}
\text{\footnotesize $300$ \normalsize}&1 & \cdot & \cdot & \cdot & \cdot & \cdot & -q & \cdot & \cdot & \frac{-q^2}{\left(q + 1\right)}\\\text{\footnotesize $210$ \normalsize}&
\cdot & 1 & -q & -q & \cdot & \frac{q^3}{\left(q + 1\right)} & -q & \frac{q^3}{\left(q + 1\right)} & -q^3 & \frac{q^3}{\left(q +\
 1\right)}\\\text{\footnotesize $201$ \normalsize}&
\cdot & \cdot & 1 & \cdot & \cdot & -q & \cdot & \frac{-q^2}{\left(q + 1\right)} & q^2 & -q\\\text{\footnotesize $120$ \normalsize}&
\cdot & \cdot & \cdot & 1 & \cdot & \frac{-q^2}{\left(q + 1\right)} & -q & -q & q^2 & \frac{q^3}{\left(q + 1\right)}\\\text{\footnotesize $111$ \normalsize}&
\cdot & \cdot & \cdot & \cdot & 1 & -q & \cdot & -q & q\left(q + 1\right) & q^2\\\text{\footnotesize $102$ \normalsize}&
\cdot & \cdot & \cdot & \cdot & \cdot & 1 & \cdot & \cdot & -q & -q\\\text{\footnotesize $030$ \normalsize}&
\cdot & \cdot & \cdot & \cdot & \cdot & \cdot & 1 & \cdot & \cdot & \frac{-q^2}{\left(q + 1\right)}\\\text{\footnotesize $021$ \normalsize}&
\cdot & \cdot & \cdot & \cdot & \cdot & \cdot & \cdot & 1 & -q & -q\\\text{\footnotesize $012$ \normalsize}&
\cdot & \cdot & \cdot & \cdot & \cdot & \cdot & \cdot & \cdot & 1 & -q\\\text{\footnotesize $003$ \normalsize}&
\cdot & \cdot & \cdot & \cdot & \cdot & \cdot & \cdot & \cdot & \cdot & 1\end{array} \, .
$$ \normalsize
\end{example}

\section{Orthogonality properties for the $q$-Key polynomials}

We show in this section that the $q$-Key polynomials $U_v$ and $\Uh_v$
are two adjoint bases with respect to a certain 
scalar product.

\subsection{A scalar product on $\Pol$}

For any Laurent series $f= \sum_{i=k}^\infty f_i x^i$, we denote
by $CT_x(f)$ the coefficient $f_0$.

Let 
$$ \Theta:=\prod_{1\leq i<j\leq n} \frac{1-x_i/x_j}{1-qx_i/x_j} \ .$$
Therefore, for any Laurent polynomial $f(x_1,\ldots,x_n)$,
the expression
$$ CT(f\, \Theta):=  CT_{x_n}\left(CT_{x_{n-1}}\left(\ldots
  \left(CT_{x_1}\left(f\, \Theta 
   \right)\right)\ldots \right)\right)
  $$
 is well defined.
 Let us use it to define a bilinear form $(\, ,\,)_q$ on $\Pol$ by
\begin{equation}
(f\, ,\, g)_q = CT\left( f\, g^\clubsuit \prod_{1\leq i<j\leq n}
  \frac{1-x_i/x_j}{1-qx_i/x_j} \right)
\end{equation} 
where $\clubsuit$ is the automorphism defined by $x_i \longmapsto
1/x_{n+1-i}$ for $1\leq i\leq n$.

Since $\Theta$ is 
invariant under $\clubsuit$, the form $(\, , \,)_q$ is
symmetrical. Under the specialization $q=0$, the previous scalar
product becomes
\begin{equation}
(f\, ,\, g) := (f\, ,\,g)\big \vert_{q=0}
 = CT\left( f\, g^\clubsuit \prod_{1\leq i<j\leq n} (1-x_i/x_j)\right).
\end{equation} 
We can also write $(f\,,\,g)_q=(f\,,\,g\Omega)$ with $\Omega=
\prod_{1\leq i<j\leq n}(1-qx_i/x_j)^{-1}$.\\

Notice that, interpreting Schur functions as characters of unitary
groups, Weyl defined the  scalar product of two  symmetric functions
$f,g$ in $n$ variables 
as the constant term of
$$ \frac{1}{n!}\, f\, g^\clubsuit\, 
               \prod_{i,j:\, i\neq j} ( 1-x_i/x_j) \, .$$ 
Essentially, Weyl takes the square of the Vandermonde,
while we are taking the quotient of the Vandermonde by the $q$-Vandermonde.

We now examine the compatibility of $\square_i$ and
$\nabla_i$ with the scalar product.

\begin{lemma}\label{adjointOperators} 
For $i$ such that $1\leq i\leq n\moins 1$, $\square_i$
(resp. $\nabla_i$) is adjoint to $\square_{n-i}$
(resp. $\nabla_{n-i}$) with respect to $(\, ,\,)_q$. 
\end{lemma}

\proof    Since $\pi_i$ (resp. $\ph_i$) is adjoint to $\pi_{n-i}$
(resp. $\ph_{n-i}$) with respect to $(\, , \,)$ (see \cite{CBMS} for
more details), we have
\begin{align*}
  ( f \square_i\,,\, g )_q &=(f\,,\,g\ \Omega\
  \pi_{n-i}(1-qx_{n-i+1}/x_{n-i}) ) \\ &=( f\,,\,g\
  \frac{(1-qx_{n-i+1}/x_{n-i})}{(1-qx_{n-i+1}/x_{n-i})}\ \Omega\
  \pi_{n-i}(1-qx_{n-i+1}/x_{n-i}) )
\end{align*}
Since the polynomial $\Omega/(1-qx_{n-i+1}/x_{n-i})$ is symmetrical in
the indeterminates $x_{n-i}$ and $x_{n-i+1}$, it commutes with the
action of $\pi_{n-i}$. Therefore 
$$(f \square_i\,,\, g)_q=(f\,,\,g\ (1-qx_{n-i+1}/x_{n-i})\ \pi_{n-i} \
\Omega)=(f\,,\, g\ \square_{n-i})_q \ . $$ 
This proves that $\square_i$ is
adjoint to $\square_{n-i}$, and, equivalently,  that $\nabla_i$ 
is adjoint to $\nabla_{n-i}$.    \QED

We shall need to characterize whether the scalar product of
two monomials vanishes or not.
Notice that, by definition, 
 $$ (x^u\, , \, x^v)=(x^{u-v\omega}\,,\,1)\, ,$$
so that one of the two monomials can be taken equal to $1$.

\begin{lemma}\label{scalarMonomials} 
 For any  $u\in \Z^n$, then 
 $(x^u, \,1)_q \neq 0$  iff $|u|=0$ and  
  $u\ge \lbrack 0,\ldots 0\rbrack$. In that case,
 $(x^u, \,1)_q=Q_u(0)$.
\end{lemma}

\proof  Let us first show  that the scalar products
$(x^u\,,\,1)_q$  satisfy the same relations
(\ref{rule1}) as  the Hall-Littlewood functions $Q_u$.\\

Let $k$ be a positive integer less than $n$. Write $x_k=y$, $x_{k+1}=z$.
Any monomial $x^v$ can be written $x^t\, y^a z^b$, with $x^t$ of degree
$0$ in $x_k, x_{k+1}$. The product 
$$
x^t (y^a z^b +y^b z^a) (z-qy)\prod_{1\leq i<j\leq n} \frac{
1-x_i/x_j}{1-qx_i/x_j}$$
is equal to 
$$  (y^a z^b +y^b z^a) (z-qy) \frac{1-y/z}{1-qy/z}\, F_1 =
  (y^a z^b +y^b z^a) (z-y) F_1 \, $$
with $F_1$ symmetrical in $y,z$. 
The constant term 
$CT_{x_{k-1}}\ldots CT_{x_1}
 (x^t (y^a z^b +y^b z^a)  F_1) =F_2$ is still symmetric
in $x_k, x_{k+1}$. Therefore
$$CT_{y} \Bigl( CT_{z}\bigl( (z-y) F_2 \bigr) \Bigr) $$
is null, and in final 
$$ CT\left (x^t (y^az^b+y^bz^a)(z-qy)\prod_{1\leq
  i<j\leq n} \frac{1-x_i/x_j}{1-qx_i/x_j}\right )=0 \, .$$
 
This relation can be rewritten
$$(y^az^{b+1}x^t,1)_q+(y^{b+1}z^{a+1}x^t,1)_q-q(y^{b+1}z^{a}x^t,1)_q-q(y^{a+1}z^bx^t,1)_q=0 \ ,
$$ which is, indeed, relation (\ref{rule1}).

On the other hand, if $u_n<0$, then there is no term of degree $0$ in $x_n$ in
 $x^u \prod_{1\leq i<j\leq n} (1-x_i/x_j) (1-qx_i/x_j)^{-1} $,
and $(x^u,1)=0$, so that rule (\ref{rule2}) is also satisfied.

In consequence, the function $u\in \Z^n \to (x^u,1)$ is determined
by the values $(x^\lambda,1)$, $\l$ partition, as the function
$u\in \Z^n \to Q_u$ is determined by its restriction to partitions.
However, for degree reasons, $(x^\lambda,1)=0$ if $\l\neq 0$. 
Since $(x^0,1)=1$, one has in final that 
$(x^u,1)= Q_u(0)$.      \QED 

\begin{example} 
For $u=[1,0,3]$ and $v=[0,1,3]$, 
$$(x^{103}\,,\,x^{013})_q=(x^{-2,-1,3}\,,\,1)_q=
Q_{-2,-1,3}(0) = q^2 (1-q) (1-q^2) \ .$$
\end{example}

\subsection{Duality between $(U_v)_{v\in\mathbb{N}^n}$ 
and $(\Uh_v)_{v\in\mathbb{N}^n}$}

Using that $\square_i$ is adjoint to $\square_{n-i}$,
we are going to prove in this section 
that $U_v$ and $\Uh_v$ are two adjoint bases of
$\Pol$ with respect to the scalar product $(\, ,\,)_q$.

We first need some technical lemmas, to allow an induction 
on the $q$-Key polynomials, starting with dominant weights.

\begin{lemma}\label{rec0}
Let $i$ be an integer such that $1\leq i\leq n\moins 1$, let
$f_1,f_2,g_1$ be three polynomials and $b$ be a constant such that
$$ f_2= f_1\, (\square_i+b) \ , \ (f_1\,,\,g_1)_q= 0
\ \ \text{and} \ \ (f_2\,,\,g_1)_q=1  \, .$$
Then the polynomial $g_2= g_1 (\nabla_{n-i} -b)$ is such that 
$$  (f_1\,,\,g_2)_q= 1\ , \ (f_2\,,\,g_2)_q=0 \, . $$
\end{lemma}

\proof  
Using that $\nabla_{n-i}$ is adjoint to $\square_i$ and that 
$\square_i \nabla_i=0$, one has
\begin{multline*} 
 (f_2\,,\,g_2)_q=(f_1\, (\square_i+b)\,,\, g_1\, (\nabla_{n-i} -b))_q=(f_1\,
(\square_i+b)\, (\nabla_{i} -b)\,,\, g_1)_q   \\
= ( f_1 ( -b(1+q)-b^2), g_1)_q = 0 \, .
\end{multline*}

Similarly, we have
\begin{align*}
  (f_1\,,\,g_2)_q&=(f_1\,,\, g_1 \,(\nabla_{n-i} -b))_q\\&=(f_1\,,
\,g_1\,(\square_{n-i}-1-q-b))_q\\ &=(f_1\,(\square_{i}+b-1-q-2b)\,,\, g_1)_q
=(f_2\,,\,g_1)_q= 1 \hfill.
\end{align*}
                                    \QED

\begin{corollary}\label{rec}
  Let $i$ be an integer such that $1\leq i\leq n\moins 1$, let $V$ be
  a vector space such that $V=V^{'}\oplus <f_1,f_2>$ with
  $f_2=f_1(\square_i+b)$ and $V^{'}$ stable under $\square_i$, and let
  $g_1$ such that
  $$(f_1\,,\,g_1)_q=0 \ \ \text{and} \ \ (f_2\,,\,g_1)_q=1 \ \ \text{and} \ \ (v\,,\,
  g_1)_q = 0, \ \forall v \in V^{'}\ . $$ Then the element
  $g_2=g_1(\nabla_{n-i}-b)$ is such that
  $$(f_2\,,\,g_2)_q=0 \ \ \text{and} \ \ (f_1\,,\,g_2)_q=1 \ \ \text{and} \ \ (v\,,\,
  g_2)_q = 0, \ \forall v \in V^{'}\ . $$
\end{corollary}

\begin{lemma}\label{lemmaZero}
Let $u$ and $\lambda$ be two dominant weights and $v$ and $\mu$ two
permutations of $u$ and $\lambda$ respectively. 
If $(x^v\,,\, x^\lambda)\not = 0$
and $(x^u\,,\,x^\mu)\not =0$ then
$$u=\lambda\quad, \quad v=\lambda\omega\quad\text{and}\quad \mu=u\omega\, . $$
\end{lemma}

\proof Using lemma \ref{scalarMonomials} , the condition
   $(x^v,x^\lambda)_q\not = 0$ and $(x^u,x^\mu)_q\not = 0$ implies two
   systems of inegalities\\\\
\begin{minipage}{8cm}
\begin{center}
$\left \lbrace \begin{array}{ccc} 
v_n            & \ge    & \lambda_1\, ,\\
v_n+v_{n-1}    & \ge    & \lambda_1+\lambda_2\, ,\\
\vdots         & \vdots & \vdots \\
v_n+\ldots+v_1 & \ge    & \lambda_1+\ldots+\lambda_n\, .
\end{array}\right .$
\end{center}
\end{minipage}
\begin{minipage}{1cm}
and
\end{minipage}
\begin{minipage}{8cm}
\begin{center}
$\left \lbrace \begin{array}{ccc} 
\mu_n                          & \ge    & u_1\, ,\\
\mu_n+\mu_{n-1}    & \ge    & u_1+u_2\, ,\\
\vdots                                     & \vdots & \vdots \\
\mu_n+\ldots+\mu_1 & \ge    & u_1+\ldots+u_n\, .
\end{array}\right .$
\end{center}
\end{minipage}\\\\
The first inequalities of the systems give $ v_n\ge \lambda_1\ge \mu_n
\ge u_1\ge v_n $. Consequently $u_1=\lambda_1=v_n=u_n$. By recursion,
using the other inequalities, one gets the lemma.\hfill \QED
\begin{corollary}
 Let $v$ be a weight and $\lambda$ a dominant weight. Then,
$$(U_v, x^\lambda)_q=\delta_{v,\lambda\omega}$$
\end{corollary}

\proof           Let $u$ be the decreasing reordering of $v$ and $\sigma$
the permutation such that $U_v=x^uY_\sigma$. As the leading term of
$U_v$ is $x^v$ and using lemma \ref{lemmaZero}, we have that
$(x^uY_\sigma\,,x^\lambda)_q\not = 0$ implies
$(x^v\,,\,x^\lambda)_q\not = 0$. By denoting $\triangle_\sigma$ the
adjoint of $Y_\sigma$ with respect to $(\, , \,)_q$, we have $
(x^uY_\sigma\,,\, x^\lambda)_q=(x^u\,,\,
x^\lambda\triangle_\sigma)_q\not=0$. As the leading term of
$x^\lambda\triangle_\sigma$ is $x^{\lambda\sigma^{'}}$, where $\lambda
\sigma^{'}$ is a permutation of $\lambda$, we obtain that $(x^u\,,\,
x^{\lambda\sigma^{'}})_q\not = 0$. Using lemma \ref{lemmaZero} we
conclude that $v=\lambda\omega$. \\ \hspace*{2cm} 
 \hfill \QED

Our main result is the following duality property between $U_v$ and $\Uh_v$.
\begin{theorem}
  The two sets of polynomials $(U_v)_{v\in\N^n}$ and
  $(\Uh_v)_{v\in\N^n}$ are two adjoint bases of $\Pol$ with respect to
  the scalar product $(\, ,\, )_q$. More precisely, they satisfy
$$ ( U_v\, ,\, \Uh_{u\omega})_q= \delta_{v,u} \ .$$
\end{theorem}

\proof  Let $\lambda$ be a dominant weight and $V$ the vector
space spanned by the $U_v$ for $v$ in $\mathcal{O}(\lambda)$. The idea
of the proof is to build by iteration the elements $(\Uh_v)_{v\in
\mathcal{O}(\lambda)}$ starting with $x^\lambda=\Uh_\lambda$. By
definition of the $q$-Key polynomials, it exists a constant $b$ such
that $ U_{\lambda \omega} = U_{\lambda\omega s_1} \
(\square_{1}+b)$. One can write the decomposition
$V=V^{'}\oplus<U_{\lambda \omega} \,,\,U_{\lambda\omega s_1}>$, with
$V'$ stable under the action of $\square_1$. Using the previous lemma,
we have that
$(U_{\lambda\omega}\,,\,x^\lambda)_q=(U_{\lambda\omega}\,,\,\Uh_\lambda)_q=1$
and
$(U_{\lambda\omega\sigma_1}\,,\,x^\lambda)_q=(U_{\lambda\omega\sigma_1}\,,\,\Uh_\lambda)_q=0$. Consequently,
by lemma \ref{rec}, the function
$x^\lambda(\nabla_{n-1}-b)=\Uh_{\lambda s_1}$ satisfy the duality
conditions
$$(U_{\lambda\omega}\,,\,\Uh_{\lambda s_1})_q=0 \quad , \quad
(U_{\lambda\omega\sigma_1}\,,\,\Uh_{\lambda s_1})_q=1 \quad \text{and}
\quad (v\,,\,\Uh_{\lambda s_1})_q=0\quad \forall v \in V^{'}\, .$$ 
By iteration, this proves that for all $u,v$, one has 
 $ ( U_v\, ,\, \Uh_{u\omega})_q= \delta_{v,u} \ .$ \hfill \QED\\

This theorem implies that the space of symmetric functions and the
linear span of dominant monomials are dual of each other, the
Hall-Littlewood functions being the basis dual to dominant monomials.

We finally mention that in the case $q=0$, one has a reproducing
kernel, as stated by the following theorem of \cite{DoubleCrystal},
 which gives another implicit
definition of the scalar product $(\, , \,)$.

\begin{theorem}
The two families of polynomials $(K_v)_{v\in\N^n}$ and
$(\Kh_v)_{v\in\N^n}$ satisfy the following Cauchy formula
\begin{equation}
\sum_{u \in \mathbb{N}^n} K_u(x) \Kh_{u\omega}(y)= \prod_{i+j
    \le n+1}\frac{1}{1-x_iy_j}\ .
\end{equation}
\end{theorem}

\end{document}